\documentclass[preprint]{elsarticle}
\usepackage{amsmath,amssymb}
\newtheorem{df}{Definition}[section]
\newtheorem{thm}[df]{Theorem}

\title{A continuous spherical wavelet transform for~$\mathcal C(\mathcal S^n)$}
\author[IIN]{I. Iglewska--Nowak\corref{cor}}
\address[IIN]{West Pomeranian University of Technology, School of Mathematics, al. Piast\'ow 17, 70--310 Szczecin, Poland, e--mail: iiglewskanowak@zut.edu.pl, phone: 0048 91 4494 826}
\cortext[cor]{Corresponding author}

\begin{document}

\begin{abstract}
In the present paper, a wavelet family over the $n$-dimensional sphere is constructed such that for each scale the wavelet is a polynomial and the inverse wavelet transform of a continuous function converges in the supremum norm.
\end{abstract}

\maketitle

\section{Introduction}

In the last decades, several wavelet constructions have been proposed for spherical signals, see e.g. \cite{ADJV02,AV,AVn,sB09,CPMDHJ05,sE11,EBCK09,mF09,mF08,FW-C,FW,GM09a,GM09b,HH09,HCM03,HI07,IIN15CWT,IIN16MR,NW}. I have shown in~\cite[Section~5]{IIN15CWT} that there are only two essentially different continuous wavelet transforms for spherical signals, namely that based on group theory~\cite{AV,AVn} and that derived from approximate identities~\cite{FGS-book,EBCK09,IIN15CWT}. In both cases, square integrable signals are considered, and the inverse transform is performed by an integral that converges in $\mathcal L^2$-norm.

The present paper is devoted to wavelet analysis of continuous functions. Since the sphere is compact, continuous functions are square integrable. The wavelet analysis I apply is that derived from approximate identities~\cite{IIN15CWT}. The novelty is that the inverse wavelet transform converges uniformly.

The paper is based on~\cite{IIN15CWT} and the same notation is used. In Section~\ref{sec:bwt} a new (with respect to that from~\cite{IIN15CWT}) definition of a wavelet family is given, and it is shown that the wavelet analysis and the wavelet synthesis can be performed in the usual way. In Section~\ref{sec:wavelet} a wavelet family is constructed such that the wavelet transform is invertible in the supremum norm. Moreover, the wavelets are polynomials.

\section{The wavelet transform}\label{sec:bwt}

This section is devoted to the bilinear wavelet transform as defined in~\cite[Section~3]{IIN15CWT}. The goal is to relax the constraints on a function family to be a wavelet, more precisely, to abandon the condition
\begin{equation}\label{eq:admbwv2}
\int_{-1}^1\left|\int_R^\infty\bigl(\overline{\Psi_\rho}\ast\Psi_\rho\bigr)(t)\,\frac{d\rho}{\rho}\right|\,\left(1-t^2\right)^{\lambda-1/2}dt\leq c
\end{equation}
uniformly in~$R$. This inequality has been verified for the Poisson wavelets in the proof of Theorem~8.2 in~\cite{IIN15PW}. Note that that proof is based on some previous results which required a deep study of the wavelet structure. I did not find in the literature any proof of~\eqref{eq:admbwv2} for the Gauss-Weierstrass wavelet, a standard example of a spherical wavelet. Here, I provide a new convergence proof for the inverse wavelet transform such that~\eqref{eq:admbwv2} can be dropped in the wavelet definition.

\begin{df}\label{def:bilinear_wavelets} Let $\alpha:\mathbb R_+\to\mathbb R_+$ be a weight function, piecewise continuous and bounded on compact sets. A family $\{\Psi_\rho\}_{\rho\in\mathbb R_+}\subseteq\mathcal L^2(\mathcal S^n)$ is called an admissible spherical wavelet  if it satisfies the following condition:
\begin{equation}\label{eq:admbwv}
\sum_{k\in\mathcal M_{n-1}(l)}\int_0^\infty\left|a_l^k(\Psi_\rho)\right|^2\cdot\alpha(\rho)\,d\rho=N(n,l)\qquad\text{for }l\in\mathbb{N}_0.
\end{equation}
\end{df}

The wavelet transform is defined by
\begin{equation}\label{eq:bwt}
\mathcal W_\Psi f(\rho,g)=\int_{\mathcal S}\overline{\Psi_\rho(g^{-1}x)}\cdot{}_sf(x)\,d\sigma(x).
\end{equation}

\begin{thm}\label{thm:inversion}If $\{\Psi_\rho\}_{\rho\in\mathbb R_+}$ is an admissible wavelet, then the wavelet transform is invertible by
\begin{equation}\label{eq:bwt_synthesis}
f(x)=\lim_{R\to 0}\int_R^{1/R}\!\!\int_{SO(n+1)}\Psi_\rho(g^{-1}x)\cdot(\mathcal W_\Psi\,f)(\rho,g)\,d\nu(g)\cdot\alpha(\rho)\,d\rho
\end{equation}
with limit in $\mathcal L^2$-sense.
\end{thm}

\begin{bfseries}Proof. \end{bfseries}Denote the inner integral in~\eqref{eq:bwt_synthesis} by $\mathfrak r_\rho(x)$. By~\eqref{eq:bwt}, it is equal to
\begin{equation}\label{eq:r_rho}
\mathfrak r_\rho(x)
   =\int_{SO(n+1)}\int_{\mathcal S^n}\Psi_\rho(g^{-1}x)\cdot\overline{\Psi_\rho(g^{-1}y)}\cdot f(y)\,d\sigma(y)\,d\nu(g).
\end{equation}
Since the functions $\Psi_\rho$ and~$f$ are square integrable, the integral is convergent and the order of integration may be changed according to Fubini's theorem,
\begin{align*}
\mathfrak r_\rho(x)&=\int_{\mathcal S^n}\int_{SO(n+1)}
   \Psi_\rho(g^{-1}x)\cdot\overline{\Psi_\rho(g^{-1}y)}\,d\nu(g)\cdot f(y)\,d\sigma(y)\\
&=\int_{\mathcal S^n}\left(\Psi_\rho\,\hat\ast\,\overline{\Psi_\rho}\right)(x\cdot y)\cdot f(y)\,d\sigma(y).
\end{align*}
According to \cite[formula~(19)]{IIN15CWT}, the zonal product of the wavelets can be expressed as
\begin{equation}\label{eq:zonal_product_Psi}
\left(\Psi_\rho\,\hat\ast\,\overline{\Psi_\rho}\right)(x\cdot y)
   =\sum_{l=0}^\infty\sum_{k\in\mathcal M_{n-1}(l)}\frac{\left|a_l^k(\Psi_\rho)\right|^2}{N(n,l)}\cdot\frac{\lambda+l}{\lambda}\,\mathcal C_l(x\cdot y).
\end{equation}
Since
\begin{equation}\label{eq:estimationC}\begin{split}
&\left|\frac{(\lambda+l)\cdot\mathcal C_l(x\cdot y)}{N(n,l)\cdot\lambda}\right|\leq\frac{1}{N(n,l)}\cdot\frac{\lambda+l}{\lambda}\cdot\mathcal C_l(1)\\
&\qquad=\frac{(n-1)!\,l!}{(n+2l-1)(n+l-2)!}\cdot\frac{n+2l-1}{n-1}\cdot\frac{(n+l-2)!}{(n-2)!\,l!}=1
\end{split}\end{equation}
(see \cite[formula~(8.937.4)]{GR} for the value of~$\mathcal C_l(1)$) and for each~$\rho$ the Fourier coefficients of~$\Psi_\rho$ are square summable, the series on the right-hand-side of~\eqref{eq:zonal_product_Psi} is uniformly bounded (with respect to $y$, the scale parameter~$\rho$ is fixed for the moment). Further, $f\in\mathcal L^2(\mathcal S^n)\subset\mathcal L(\mathcal S^n)$. Thus, the order of integration over~$\mathcal S^n$ and summation with respect to~$l$ in
$$
\mathfrak r_\rho(x)=\int_{\mathcal S^n}\,\sum_{l=0}^\infty\sum_{k\in\mathcal M_{n-1}(l)}\frac{\left|a_l^k(\Psi_\rho)\right|^2}{N(n,l)}
   \cdot\frac{\lambda+l}{\lambda}\,\mathcal C_l(x\cdot y)\cdot f(y)\,d\sigma(y)
$$
can be changed. By the reproducing property of the kernels $\frac{\lambda+l}{\lambda}\,\mathcal C_l$, $l\in\mathbb N_0$, we obtain
\begin{align}
\mathfrak r_\rho(x)&=\sum_{l=0}^\infty\sum_{k\in\mathcal M_{n-1}(l)}\frac{\left|a_l^k(\Psi_\rho)\right|^2}{N(n,l)}\cdot f_l(x)\notag\\
&=\sum_{l=0}^\infty
   \left(\sum_{k\in\mathcal M_{n-1}(l)}\frac{\left|a_l^k(\Psi_\rho)\right|^2}{N(n,l)}\cdot\sum_{m\in\mathcal M_{n-1}(l)}a_l^m(f)\cdot Y_l^m(x)\right).\label{eq:r_rho_series}
\end{align}
Consider the coefficients
$$
(\widehat{\mathfrak r_\rho})_l^m:=\sum_{k\in\mathcal M_{n-1}(l)}\left|a_l^k(\Psi_\rho)\right|^2\cdot\frac{a_l^m(f)}{N(n,l)}.
$$
Since $f\in\mathcal L^2(\mathcal S^n)$, the Fourier coefficients $a_l^m(f)$ are bounded. Further, $N(n,l)$ is the number of summands with index~$l$ and, consequently, $\sum_m\frac{a_l^m(f)}{N(n,l)}$ is bounded. On the other hand, $\sum_k\left|a_l^k(\Psi_\rho)\right|^2$ is summable with respect to~$l$ by the square summability of the Fourier coefficients. Thus, the coefficients $(\widehat{\mathfrak r_\rho})_l^m$ are summable. Consequently, they are square summable and $\mathfrak r\in\mathcal L^2(\mathcal S^n)$. This justifies the usage of notation~$(\widehat{\mathfrak r_\rho})_l^k$. By the Plancherel theorem, for any~$R$, the $\mathcal L^2$-norm $L(R)$ of
$$
f-\int_R^{1/R}\mathfrak r_\rho\cdot\alpha(\rho)\,d\rho
$$
satisfies
\begin{equation}\label{eq:L(R)}
[L(R)]^2=\sum_l\sum_{m}
   \left[1-\int_R^{1/R}\sum_{k}\frac{\left|a_l^k(\Psi_\rho)\right|^2}{N(n,l)}\cdot\alpha(\rho)\,d\rho\right]^2\cdot\left[a_l^m(f)\right]^2.
\end{equation}
Now, by~\eqref{eq:admbwv}, for each $(l,k)$ the difference in brackets on the right-hand-side of~\eqref{eq:L(R)} is bounded. Thus, by the square summability of~$\left(a_l^m(f)\right)_{l,k}$, the series in~\eqref{eq:L(R)} is uniformly bounded. Consequently,
\begin{align*}
&\lim_{R\to0}[L(R)]^2\\
&=\sum_l\sum_{m}\lim_{R\to0}
   \left[1-\int_R^{1/R}\frac{1}{N(n,l)}\sum_{k}\left|a_l^k(\Psi_\rho)\right|^2\cdot\alpha(\rho)\,d\rho\right]^2\cdot\left[a_l^m(f)\right]^2=0,
\end{align*}
and the assertion follows.\hfill$\Box$\\[-1em]

Similarly, the isometry condition can be proven for wavelets defined by~\eqref{eq:admbwv}.

\begin{thm}\label{thm:isometry} Let $\{\Psi_\rho\}_{\rho\in\mathbb R_+}$ be an admissible spherical wavelet  and $f,g\in\mathcal L^2(\mathcal S^n)$. Then,
$$
\left<\mathcal W_\Psi\,f,\mathcal W_\Psi\,g\right>=\left<f,g\right>
$$
where the scalar product in the wavelet phase space is given by
$$
\left<F,G\right>_{\mathcal L^2(\mathbb R_+\times SO(n+1))}=\int_0^\infty\!\!\int_{SO(n+1)}\overline{F(\rho,g)}\,G(\rho,g)\,d\nu(g)\,\alpha(\rho)\,d\rho.
$$
\end{thm}

The proof is analogous to that of Theorem~\ref{thm:inversion}.

\begin{bfseries}Proof.\end{bfseries} For a fixed~$\rho$, Denote by~$\mathfrak r_\rho$ the integral
\begin{align*}
\mathfrak r_\rho&=\int_{SO(n+1)}\overline{(\mathcal W_\Psi\,f)(\rho,g)}\cdot(\mathcal W_\Psi\,g)(\rho,g)\,d\nu(g)\\
&=\int_{SO(n+1)}\int_{\mathcal S^n}\Psi_\rho(g^{-1}x)\,\overline{f(x)}\,d\sigma(x)
   \cdot\int_{\mathcal S^n}\overline{\Psi_\rho(g^{-1}y)}\,g(y)\cdot d\sigma(y)\,d\nu(g).
\end{align*}
A change of the integration order yields
\begin{align*}
\mathfrak r_\rho=\int_{\mathcal S^n}\int_{\mathcal S^n}&
   \underbrace{\int_{SO(n+1)}\Psi_\rho(g^{-1}x)\cdot\overline{\Psi_\rho(g^{-1}y)}\cdot d\nu(g)}
   _{\Psi_\rho\,\hat\ast\,\overline{\Psi_\rho}(x,y)}\\
&\cdot\alpha(\rho)\,d\rho\cdot\overline{f(x)}\,d\sigma(x)\cdot g(y)\,d\sigma(y).
\end{align*}
According to~\eqref{eq:zonal_product_Psi}, this expression is equal to
$$
\mathfrak r_\rho=\int_{\mathcal S^n}\int_{\mathcal S^n}\sum_{l=0}^\infty\sum_{k\in\mathcal M_{n-1}(l)}\frac{\left|a_l^k(\Psi_\rho)\right|^2}{N(n,l)}
   \cdot\frac{\lambda+l}{\lambda}\,\mathcal C_l(x\cdot y)\cdot\overline{f(x)}\,d\sigma(x)\cdot g(y)\,d\sigma(y).
$$
The series is uniformly bounded with respect to~$x$. Thus, the order of integration and summation may be changed,
\begin{align*}
\mathfrak r_\rho&=\int_{\mathcal S^n}\sum_{l,k}\frac{\left|a_l^k(\Psi_\rho)\right|^2}{N(n,l)}
   \cdot\int_{\mathcal S^n}\frac{\lambda+l}{\lambda}\,\mathcal C_l(x\cdot y)\cdot\overline{f(x)}\,d\sigma(x)\cdot g(y)\,d\sigma(y)\\
&=\int_{\mathcal S^n}\sum_{l,k}\frac{\left|a_l^k(\Psi_\rho)\right|^2}{N(n,l)}\cdot\overline{f_l(y)}\cdot g(y)\,d\sigma(y).
\end{align*}
Since by the H\"older inequality and~\eqref{eq:estimationC},
$$
\left|\frac{f_l(y)}{N(n,l)}\right|\leq\left\|\frac{\lambda+l}{\lambda\cdot N(n,l)}\,\mathcal C_l(\circ,y)\right\|_{\mathcal L^2}\cdot\|f\|_{\mathcal L^2}\leq C,
$$
and $\Psi_\rho\in\mathcal L^2(\mathcal S^n)$, the series is absolutely convergent, and the order of summation and integration may be changed,
$$
\mathfrak r_\rho=\sum_{l,k}\frac{\left|a_l^k(\Psi_\rho)\right|^2}{N(n,l)}
   \cdot\int_{\mathcal S^n}\overline{f_l(y)}\cdot g(y)\,d\sigma(y).
$$
By the orthogonality of the spherical harmonics,
$$
\int_{\mathcal S^n}\overline{f_l(y)}\cdot g(y)\,d\sigma(y)=\sum_{m\in\mathcal M_{n-1}(l)}\overline{a_l^m(f)}\cdot a_l^m(g).
$$
The assertion follows by the square summability of the Fourier coefficients $a_l^m(f)$, $a_l^m(g)$ and condition~\eqref{eq:admbwv}, with the same arguments as in the proof of Theorem~\ref{thm:inversion}.\hfill $\Box$

\section{The wavelet construction}\label{sec:wavelet}

The goal of this section is to define a wavelet family such that the integral~\eqref{eq:bwt_synthesis} converges not only in $\mathcal L^2$-norm, but also uniformly. Moreover, for each scale~$\rho$ the wavelet is a polynomial. The key to the construction is the Poisson integral formula \cite[Theorem~3.4.1]{FGS-book} and its implications \cite[Theorem~3.4.2 and Corollary~3.4.3]{FGS-book}. These statements from~\cite{FGS-book} can be straightforward generalized to the case of the $n$-dimensional sphere (with the $n$-dimensional Poisson kernel, the Gegenbauer polynomials in place of the Legendre polynomials, and spherical harmonics replaced by the hyperspherical ones) and we omit their proofs here.

\begin{thm}\label{thm:wavelet_for_C}Let $\alpha(\rho)=\frac{1}{\rho}$ and let $\{\Psi_\rho\}_{\rho\in\mathbb R_+}$ be a wavelet family defined by
\begin{equation}\label{eq:wavelet_coeffs}\begin{split}
a_l^0(\Psi_\rho)&=\sqrt{\rho\cdot N(n,l)\cdot\gamma_l(\rho)},\\
a_l^k(\Psi_\rho)&=0\qquad\text{for }k\ne(0,0,\dots,0).
\end{split}\end{equation}
with~$\gamma_l$, $l\in\mathbb N_0$, given by
\begin{align*}
\gamma_0(\rho)&=\chi_{[1,2]}(\rho),\\
\gamma_l(\rho)
   &=\begin{cases}
      0,&\rho\in\left[\frac{1}{l},\infty\right),\\
      \frac{l^{l+1}}{(l+1)^{l-1}},&\rho\in\left[\frac{1}{l+1},\frac{1}{l}\right),\\
      l(1-\rho)^{l-1},&\rho\in\left(0,\frac{1}{l+1}\right)
   \end{cases}
   \qquad\text{for }l>0.
\end{align*}
Then for each $f\in\mathcal C(\mathcal S^n)$ the right-hand-side of~\eqref{eq:bwt_synthesis} converges uniformly to~$f$.
\end{thm}

\begin{bfseries}Remark. \end{bfseries}Note that the wavelet family defined in Theorem~\ref{thm:wavelet_for_C} is a zonal one such that the wavelet analysis and the wavelet synthesis can be simplified according to the formulae given in \cite[Subsection~3.2]{IIN15CWT}.

\begin{bfseries}Proof. \end{bfseries}The functions~$\gamma_l$, $l\in\mathbb N$, satisfy
$$
\int_R^\infty\gamma_l(\rho)\,d\rho=\begin{cases}(1-R)^l, &R\in\left[0,\frac{1}{l+1}\right],\\0, &R\in\left[\frac{1}{l},\infty\right)\end{cases}
$$
and for $R\in\left[\frac{1}{l+1},\frac{1}{l}\right]$ the integral $\int_R^\infty\gamma_l(\rho)\,d\rho$ has values between~$0$ and~$(1-R)^l$. Further,
$$
\int_R^\infty\gamma_0(\rho)\,d\rho=1\qquad\text{for }R\in[0,1].
$$
Thus, $a_l^k(\Psi_\rho)$ satisfy~\eqref{eq:admbwv}, i.e., $\{\Psi_\rho\}$ is a wavelet. Consider~$\mathfrak r_\rho$ given by~\eqref{eq:r_rho}. According to~\eqref{eq:r_rho_series} it is equal to
$$
\mathfrak r_\rho(x)=\sum_{l=0}^\infty\rho\cdot\gamma_l(\rho)\sum_{m\in\mathcal M_{n-1}(l)}a_l^m(f)\cdot Y_l^m(x).
$$
Thus, for $R\in\left[\frac{1}{L+1},\frac{1}{L}\right]$, $L\in\mathbb N$,
\begin{align*}
\int_R^{1/R}\mathfrak r_\rho(x)\,\frac{d\rho}{\rho}&=\sum_{l=0}^{L-1}(1-R)^lf_l(x)+\alpha\,(1-R)^Lf_L(x).
\end{align*}
with $\alpha\in[0,1]$. The assertion follows similarly as in the proof of \cite[Corollary~3.4.3]{FGS-book}.\hfill$\Box$


\begin{thebibliography}{99}
\bibitem{ADJV02} J.-P. Antoine, L. Demanet, L. Jacques, and P. Vandergheynst, \emph{Wavelets on the sphere: implementation and approximations}, Appl. Comput. Harm. Anal. 13 (2002), 177--200.
\bibitem{AV} J.-P. Antoine and P. Vandergheynst, \emph{Wavelets on the 2-sphere: a group-theoretical approach}, Appl. Comput. Harmon. Anal. 7 (1999), No.~3, 262--291.
\bibitem{AVn} J.-P. Antoine and P. Vandergheynst, \emph{Wavelets on the n-sphere and related manifolds}, J. Math. Phys. 39 (1998), No.~8, 3987--4008.
\bibitem{sB09} S.~Bernstein, \emph{Spherical singular integrals, monogenic kernels and wavelets on the three--dimensional sphere}, Adv. Appl. Clifford Algebr. 19 (2009), No.~2, 173--189.
\bibitem{CPMDHJ05} A. Chambodut, I. Panet, M. Mandea, M. Diament, M. Holschneider, and O. Jamet, \emph{Wavelet frames: an alternative to spherical harmonic representation of potential fields}, Geophys. J. Int. 163 (2005), 875--899.
\bibitem{sE11} S. Ebert, \emph{Wavelets on Lie groups and homogeneous spaces}, PhD--thesis, Freiberg 2011.
\bibitem{EBCK09} S. Ebert, S. Bernstein, P. Cerejeiras, and U. K\'ahler, \emph{Nonzonal wavelets on~$\mathcal S^N$}, 18$^\text{th}$ International Conference on the Application of Computer Science and Mathematics in Architecture and Civil Engineering, Weimar 2009.
\bibitem{mF09} M. Ferreira, \emph{Spherical continuous wavelet transforms arising from sections of the Lorentz group}, Appl. Comput. Harm. Anal. 26 (2009), No.~2, 212--229.
\bibitem{mF08} M. Ferreira, \emph{Spherical wavelet transform}, Adv. Appl. Clifford Alg. 18 (2008), 611--619.
\bibitem{FGS-book} W. Freeden, T. Gervens, and M. Schreiner, \emph{Constructive approximation on the sphere. With applications to geomathematics}, Clarendon Press, New York 1998.\bibitem{FW-C} W.~Freeden and U.~Windheuser, \emph{Combined spherical harmonic and wavelet expansion -- a future concept in {E}arth's gravitational determination}, Appl. Comput. Harmon. Anal. 4 (1997), No. 1, 1--37.
\bibitem{FW} W.~Freeden and U.~Windheuser, \emph{Spherical wavelet transform and its discretization}, Adv. Comput. Math. 5 (1996), No.~1, 51--94.
\bibitem{GM09a}  D. Geller and A. Mayeli, \emph{Continuous wavelets on compact manifolds}, Math. Z. 262 (2009), no. 4, 895--927.
\bibitem{GM09b}  D. Geller and A. Mayeli, \emph{Nearly tight frames and space--frequency analysis on compact manifolds}, Math. Z. 263 (2009), no. 2, 235--264.
\bibitem{HH09} M.~Hayn and M.~Holschneider, \emph{Directional spherical multipole wavelets}, J. Math. Phys. 50 (2009), No.~7, 073512, 11 p.
\bibitem{mH96} M.~Holschneider, \emph{Continuous wavelet transforms on the sphere}, J. Math. Phys. 37 (1996), No.~8, 4156--4165.
\bibitem{HCM03} M.~Holschneider, A.~Chambodut, and M.~Mandea, \emph{From global to regional analysis of the magnetic field on the sphere using wavelet frames}, Phys. Earth Planet. Inter. 135 (2003), 107--123.
\bibitem{HI07} M. Holschneider and I. Iglewska--Nowak, \emph{Poisson wavelets on the sphere}, J. Four. Anal. Appl. 13 (2007), 405--419.
\bibitem{IIN15CWT} I. Iglewska--Nowak, \emph{Continuous wavelet transforms on $n$-dimensional spheres}, Appl. Comput. Harmon. Anal. 39 (2015), no. 2, 248--276.
\bibitem{IIN16MR} I. Iglewska--Nowak, \emph{Multiresolution on $n$-dimensional spheres}, Kyushu J. Math., \emph{in press}.
\bibitem{IIN15PW} I. Iglewska--Nowak, \emph{Poisson wavelets on $n$-dimensional spheres}, J. Fourier Anal. Appl. 21 (2015), no. 1, 206--227.
\bibitem{NW} F.J.~Narcowich and J.D.~Ward, \emph{Nonstationary wavelets on the $m$-sphere for scattered data}, Appl. Comput. Harmon. Anal. 3 (1996), No.~4, 324--336.
\end{thebibliography}
\end{document}